\documentclass[preprint,12pt]{elsarticle}

\usepackage{amssymb}
\usepackage{amsmath}
\usepackage{amsthm}
\usepackage{mathrsfs}
\usepackage[colorlinks]{hyperref}
\AtBeginDocument{%
  \hypersetup{%
    linkcolor=blue,%
    citecolor=green,%
  }%
}
\usepackage{cleveref}

\usepackage[displaymath,mathlines,pagewise]{lineno}
\newcommand*\patchAmsMathEnvironmentForLineno[1]{%
  \expandafter\let\csname old#1\expandafter\endcsname\csname #1\endcsname
  \expandafter\let\csname oldend#1\expandafter\endcsname\csname end#1\endcsname
  \renewenvironment{#1}%
     {\linenomath\csname old#1\endcsname}%
     {\csname oldend#1\endcsname\endlinenomath}}%
\newcommand*\patchBothAmsMathEnvironmentsForLineno[1]{%
  \patchAmsMathEnvironmentForLineno{#1}%
  \patchAmsMathEnvironmentForLineno{#1*}}%
\AtBeginDocument{%
\patchBothAmsMathEnvironmentsForLineno{equation}%
\patchBothAmsMathEnvironmentsForLineno{align}%
\patchBothAmsMathEnvironmentsForLineno{flalign}%
\patchBothAmsMathEnvironmentsForLineno{alignat}%
\patchBothAmsMathEnvironmentsForLineno{gather}%
\patchBothAmsMathEnvironmentsForLineno{multline}%
}

\crefname{equation}{}{}

\newtheorem{theorem}{Theorem}[section]
\newtheorem{lemma}[theorem]{Lemma}

\newtheorem{corollary}[theorem]{Corollary}

\theoremstyle{definition}
\newtheorem{definition}[theorem]{Definition}

\theoremstyle{remark}
\newtheorem{remark}[theorem]{Remark}

\numberwithin{equation}{section}

\journal{~}

\begin{document}

\begin{frontmatter}



\title{Pointwise Boundary Differentiability for Fully Nonlinear Elliptic Equations\tnoteref{t1}}


\author[rvt]{Duan Wu}
\ead{wuduan@mail.nwpu.edu.cn}
\author[rvt1]{Yuanyuan Lian\corref{cor1}}
\ead{lianyuanyuan@sjtu.edu.cn; lianyuanyuan.hthk@gmail.com}
\author[rvt1]{Kai Zhang}
\ead{zhangkaizfz@gmail.com}
\tnotetext[t1]{This research is supported by the National Natural Science Foundation of China (Grant No. 11701454) and the China Postdoctoral Science Foundation (Grant No. 2021M692086).}

\cortext[cor1]{Corresponding author.}

\address[rvt]{School of Mathematics and Statistics, Northwestern Polytechnical University, Xi'an, Shaanxi, 710129, PR China}
\address[rvt1]{School of Mathematical Sciences, Shanghai Jiao Tong University, Shanghai, China}

\begin{abstract}
In this paper, we prove the pointwise boundary differentiability for viscosity solutions of fully nonlinear elliptic equations. This generalizes the previous related results for linear equations and the geometrical conditions in this paper are pointwise and more general. Moreover, our proofs are relatively simple.
\end{abstract}


\end{frontmatter}


\section{Introduction}\label{S1}
We study the pointwise boundary differentiability for viscosity solutions of
\begin{equation}\label{1.1}
\left\{\begin{aligned}
&u\in S^*(\lambda,\Lambda,f)&& ~~\mbox{in}~~\Omega;\\
&u=g&& ~~\mbox{on}~~\partial \Omega,
\end{aligned}\right.
\end{equation}
where $\Omega\subset {R^n}$ is a bounded domain (i.e. connected open set); $S^*\left( {\lambda,\Lambda,f} \right)$ is the Pucci's class (see \Cref{SP} for details). In this paper, solutions always indicate viscosity solutions.

The boundary differentiability has been studied extensively for linear uniformly elliptic equations. Li and Wang \cite{MR2264018,MR2494685} proved the boundary differentiability on convex domains. Later, Li and Zhang \cite{MR3029135} extended the boundary differentiability to more general domains, called $\gamma$-convex domains (see also \Cref{r-1.2} for the definition). Both ``convexity'' and ``$\gamma$-convexity'' are global geometrical conditions. The pointwise regularity is more attractive than the local and global regularity because it shows clearly the behavior of solutions influenced by the coefficients and the prescribed data. Moreover, the assumptions for pointwise regularity are usually weaker.

Huang, Li and Wang \cite{MR3167627} gave a pointwise geometrical condition under which they obtained the boundary pointwise differentiability. Huang, Li and Zhang \cite{MR3556828} also considered the pointwise differentiability and constructed a counterexample, demonstrating the necessity of the exterior $C^{1,\mathrm{Dini}}$ condition (see \Cref{de12}). However, all above results are obtained for linear equations. In this paper, we prove the boundary pointwise differentiability for viscosity solutions of $\cref{1.1}$ under more general pointwise geometrical conditions.

The following two theorems are our main results, which correspond to cases of the corner point and the flat point in previous papers (see \cite{MR3556828}, \cite{MR2264018,MR2494685,MR3029135}).
\begin{theorem}\label{th11}
Let $u$ be a viscosity solution of
\begin{equation*}
\left\{\begin{aligned}
&u\in S^*(\lambda,\Lambda,f)&& ~~\mbox{in}~~\Omega\cap B_1;\\
&u=g&& ~~\mbox{on}~~\partial \Omega \cap B_1,
\end{aligned}\right.
\end{equation*}
where $f\in L^{n}(\Omega\cap B_1)$ and $g$ is differentiable at $0\in \partial \Omega$. Suppose that $\Omega$ satisfies the exterior $C^1$ condition at $0$ and there exists a cone $K$ with $0$ as the vertex such that $K\cap B_1\subset \Omega^c \cap R^n_+$.

Then $u$ is differentiable at $0$ and $\nabla u(0)=\nabla g(0)$.
\end{theorem}

\begin{remark}\label{re1.6}
An open convex set $K\subset R^n$ is called a cone with $0$ as the vertex if $rx\in K $ for any $x\in K$ and $r>0$.
\end{remark}

\begin{remark}\label{re1.4}
For $u\in S^*(\lambda,\Lambda,f)$ and $f\in L^n$, only the $C^{\alpha}$ ($0<\alpha<1$) regularity can be obtained in the interior (see \cite[Chapter 4]{MR1351007} and \cite{MR882838}). In contrast, \Cref{th11} shows that one can obtain higher regularity on the boundary.
\end{remark}

\begin{theorem}\label{th12}
Let $u$ be a viscosity solution of
\begin{equation*}
\left\{\begin{aligned}
&u\in S^*(\lambda,\Lambda,f)&& ~~\mbox{in}~~\Omega\cap B_1;\\
&u=g&& ~~\mbox{on}~~\partial \Omega \cap B_1,
\end{aligned}\right.
\end{equation*}
where $f\in C^{-1,\mathrm{Dini}}(0)$ and $g\in C^{1,\mathrm{Dini}}(0)$ where $0\in \partial \Omega$. Suppose that $\Omega$ satisfies the exterior $C^{1,\mathrm{Dini}}$ condition and the interior ${C^1}$ condition at $0$.

Then $u$ is differentiable at $0$.
\end{theorem}

\begin{remark}\label{re12}
The geometrical conditions in \Cref{th11} and \Cref{th12} are more general than the ones used in previous studies. Indeed, if $\Omega$ is a $\gamma$-convex domain (including convex domains as special cases, see \Cref{r-1.2}) and $0\in \partial \Omega$ is a corner point, $\partial\Omega$ satisfies the conditions in \Cref{th11} (see \cite[Lemma 2.6]{MR3029135}); if $0$ is a flat point, $\partial\Omega$ satisfies the conditions in \Cref{th12} (see \cite[Lemma 2.7]{MR3029135}).

Huang, Li and Zhang \cite{MR3556828} proposed a pointwise geometrical condition for the pointwise boundary differentiability, which is a special case of ours. In a word, even for linear equations, all previous results with respect to boundary differentiability are covered by this paper.
\end{remark}

\begin{remark}\label{re-mod1}
As for the gradient regularity at the boundary for fully nonlinear elliptic equations, Silvestre and Sirakov \cite{MR3246039} proved the boundary $C^{1,\alpha}$ regularity. That paper requires $\partial \Omega\in C^2$ since the method of flattening the boundary is used. This result was improved by Lian and Zhang \cite{MR4088470} based on another technique without flattening the boundary. In the case of boundary differentiability, since the geometrical conditions from the interior and the exterior are not the same, above techniques are not applicable. We have to make the estimate more carefully in the proof.
\end{remark}

As a corollary, we have
\begin{corollary}\label{c-1}
Let $u$ be a viscosity solution of
\begin{equation*}
\left\{\begin{aligned}
&u\in S^*(\lambda,\Lambda,f)&& ~~\mbox{in}~~\Omega;\\
&u=g&& ~~\mbox{on}~~\partial \Omega,
\end{aligned}\right.
\end{equation*}
where $f\in C^{-1,\mathrm{Dini}}(\bar{\Omega})$ and $g\in C^{1,\mathrm{Dini}}(\bar{\Omega})$. Suppose that $\Omega$ is a $\gamma$-convex domain.

Then $u$ is differentiable on the boundary $\partial \Omega$. That is, for any $x_0\in \partial \Omega$, there exist $l\in R^n$ and a modulus of continuity $\omega_{x_0}$ such that
\begin{equation*}
|u(x)-u(x_0)-l\cdot (x-x_0)|\leq |x-x_0|\omega_{x_0}(|x-x_0|).
\end{equation*}
\end{corollary}

The main idea of this paper is the following. For understanding the boundary differentiability at $0\in\partial \Omega$ influenced by the geometrical property of $\partial \Omega$ clearly, we assume that $f,g\equiv 0$. Roughly speaking, if $\Omega$ occupies a smaller portion in a ball $B_r$ (e.g. $|\Omega\cap B_r|/|B_r|$ is smaller), the regularity of $u$ near $0$ is higher. In particular, if $\Omega$ is contained in a half ball (i.e. $\Omega\cap B_r\subset B_r^+$), it is well known that $u$ is Lipschitz continuous at $0$, which can be shown by constructing a proper barrier.

The boundary differentiability can be regarded as the enhancement of the boundary Lipschitz regularity in some sense. The treatments are different under different geometrical conditions. In \Cref{th11}, since there exists a cone $K\subset \Omega^c$ in the half ball $B_1^{+}$, this leads to a higher regularity (than the Lipschitz regularity) of $u$, namely $\nabla u(0)=0$. If $\partial \Omega$ is not curved toward the exterior too much, $\nabla u(0)=0$ will be retained.

In \Cref{th12}, the exterior $C^{1,\mathrm{Dini}}$ condition is an appropriate perturbation of the half ball and guarantees the boundary Lipschitz regularity, i.e., the boundedness of $\nabla u(0)$. In the concrete proof (see \Cref{th33}), the graph of $u$ is controlled by two hyperplanes. In addition, the interior $C^{1}$ condition will result in the convergence of the two hyperplanes to one hyperplane, which implies the boundary differentiability of $u$ at $0$.

At the end of this section, we make some comments on the boundary differentiability, which is a critical case in the theory of boundary regularity. It is well known that the boundary regularity depends heavily on geometrical properties of the domain. For the continuity of solutions up to the boundary, a geometrical condition on the boundary from the exterior is enough. Famous examples are the Wiener criterion for the boundary continuity \cite{Wiener-1924}, the exterior cone condition for the boundary H\"{o}lder continuity (see \cite[Problem 2.12, Theorem 8.27, Corollary 9.28.]{MR1814364}) and the exterior sphere condition for the boundary Lipschitz continuity (see \cite[Chapter 2.8]{MR1814364} and \cite[Lemma 1.2]{Safonov2008}) etc. However, for the boundary differentiability, only a condition from the exterior is not enough (see the counterexample in \cite{MR3556828}), which is different from the boundary continuity.

In addition, for higher order regularity, we need the same geometrical conditions on the boundary from both interior and exterior sides. For instance, the $C^{1,\alpha}$ smoothness of the boundary implies the boundary $C^{1,\alpha}$ regularity (see \cite[Theorem 8.33]{MR1814364}, \cite{MR4088470} and \cite{MR2853528}) and the $C^{2,\alpha}$ smoothness of the boundary implies the boundary $C^{2,\alpha}$ regularity (see \cite[Theorem 6.6]{MR1814364} and \cite{MR4088470}). Whereas, for the boundary differentiability, we don't need the same geometrical conditions from both sides.

In the next section, we gather some notations and preliminary results. In \Cref{Sth11}, we will give the proof of \Cref{th11} and \Cref{th12} will be proved in \Cref{Sth12}.

\section{Preliminaries}\label{SP}

\subsection{Notations and notions}\label{SP-1}
In this paper, we use standard notations. Let $\{e_i\}^{n}_{i=1}$ denote the standard basis of $R^n$ and $R^n_+=\{x\in R^n\big|x_n>0\}$. We may write $x=(x_1,...,x_n)=(x',x_n)$ for $x\in R^{n}$.
The $B_r(x_0)$ (or $B(x_0,r)$) is the open ball with radius $r$ and center $x_0$ and $B_r^+(x_0)=B_r(x_0)\cap R^n_+$. For simplicity, we denote $B_r=B_r(0)$ and $B_r^+=B^+_r(0)$. We also use the ball in $R^{n-1}$, i.e., $T_r(x_0)=B_r(x_0)\cap \left\{x_n=0\right\}$. Similarly, we may write $T_r=T_r(0)$ for short. In addition, we introduce two new notations for convenience: $\Omega_r=\Omega\cap B_r$ and $(\partial\Omega)_r=\partial\Omega\cap B_r$.

Next, we collect some notions with respect to fully nonlinear elliptic equations and viscosity solutions (see \cite{MR1351007,MR1376656,MR1118699} for details). We call $F:S^n\rightarrow R$ a fully nonlinear uniformly elliptic operator with ellipticity constants $0<\lambda\leq \Lambda$ if
\begin{equation*}
    \lambda\|N\|\leq F(M+N)-F(M)\leq \Lambda\|N\|,~~~~\forall~M,N\in S^n,~N\geq 0,
\end{equation*}
where $S^n$ denotes the set of $n\times n$ symmetric matrices; $\|N\|$ is the spectral radius of $N$ and $N\geq 0$ means the nonnegativeness.

Throughout this article, a constant is called universal if it depends only on the dimension $n$ and the ellipticity constants $\lambda$ and $\Lambda$. Unless stated otherwise, the letters $C$ and $c$ always denote positive universal constants and $0<c<1$.

\begin{definition}\label{d-viscoF}
Let $u\in C(\Omega)$ and $f\in L^{n}(\Omega)$. We say that $u$ is an $L^n$-viscosity subsolution (resp., supersolution) of
\begin{equation}\label{FNE}
F(D^2u)=f ~~~~\mbox{in}~\Omega,
\end{equation}
if
\begin{equation*}
  \begin{aligned}
    \mathrm{ess}~\underset{y\to x}{\lim\sup}\left(F(D^2\varphi(y))-f(y)\right)\geq 0\\
    \left(\mathrm{resp.},~\mathrm{ess}~\underset{y\to x}{\lim\inf}\left(F(D^2\varphi(y))-f(y)\right)\leq 0\right)
  \end{aligned}
\end{equation*}
provided that for $\varphi\in W^{2,n}(\Omega)$, $u-\varphi$ attains its local maximum (resp., minimum) at $x\in\Omega$.

We call $u$ an $L^n$-viscosity solution of \cref{FNE} if it is both an $L^n$-viscosity subsolution and supersolution of \cref{FNE}.
\end{definition}

\begin{remark}\label{r-14}
In fact, one can define $L^p$-viscosity solution for any $p>n/2$ \cite{MR1376656}. In this paper, we only deal with the $L^n$-viscosity solution and a viscosity solution always means the $L^n$-viscosity solution. If $\varphi\in W^{2,n}(\Omega)$ is replaced by $\varphi\in C^2(\Omega)$, we arrive at the definition of $C$-viscosity solution, which is adopted in \cite{MR1351007}. If all functions are continuous in their variables, these two notions of viscosity solutions are equivalent (see \cite[Proposition 2.9]{MR1376656}).
\end{remark}

The following notions of the Pucci's extremal operators and the Pucci's class are our main concerns for the boundary differentiability.
\begin{definition}\label{d-Sf}
For $M\in S^n$, denote its eigenvalues by $\lambda_i$ ($1\leq i\leq n$) and define the Pucci's extremal operators:
\begin{equation*}
\begin{aligned}
  &\mathcal{M}^+(M,\lambda,\Lambda)=\Lambda\left(\sum_{\lambda_i>0}\lambda_i\right)
+\lambda\left(\sum_{\lambda_i<0}\lambda_i\right), \\
  &\mathcal{M}^-(M,\lambda,\Lambda)=\lambda\left(\sum_{\lambda_i>0}\lambda_i\right)
+\Lambda\left(\sum_{\lambda_i<0}\lambda_i\right),
\end{aligned}
\end{equation*}
which are two typical fully nonlinear uniformly elliptic operators.

Then we can define the Pucci's class as follows. We say that $u\in \underline{S}(\lambda,\Lambda,f)$ if $u$ is an $L^n$-viscosity subsolution of
\begin{equation*}
  \mathcal{M}^+(D^2u,\lambda,\Lambda)= f.
\end{equation*}
Similarly, we denote $u\in \bar{S}(\lambda,\Lambda,f)$ if $u$ is an $L^n$-viscosity supersolution of
\begin{equation*}
  \mathcal{M}^-(D^2u,\lambda,\Lambda)= f.
\end{equation*}
We also define
\begin{equation}\label{SC2Sf}
\begin{aligned}
  &S(\lambda,\Lambda,f)=\underline{S}(\lambda,\Lambda,f)\cap \bar{S}(\lambda,\Lambda,f),\\
  &S^*(\lambda,\Lambda,f)=\underline{S}(\lambda,\Lambda,-|f|)\cap \bar{S}(\lambda,\Lambda,|f|).
\end{aligned}
\end{equation}

Clearly, $S(\lambda,\Lambda,f)\subset S^*(\lambda,\Lambda,f)$ and $S(\lambda,\Lambda,0)=S^*(\lambda,\Lambda,0)$. One feature of the Pucci's class is that for any viscosity solution $u$ of \cref{FNE}, $u\in S(\lambda,\Lambda,f)$.
\end{definition}

Now we present some geometrical conditions on the domain, under which we will prove the boundary differentiability. Recall that $\omega:[0,+\infty)\rightarrow [0,+\infty)$ is called a modulus of continuity if $\omega$ is a nonnegative nondecreasing function satisfying $\omega(r)\rightarrow 0$ as $r\rightarrow 0$. Moreover, a modulus of continuity $\omega$ is called a Dini function if it satisfies the following Dini condition for some $r_0>0$
\begin{equation}\label{deeq11}
    \int_{0}^{r_0}\frac{\omega(r)}{r}dr<\infty.
\end{equation}

\begin{definition}[\textbf{exterior $C^{1}$ ($C^{1,\mathrm{Dini}}$) condition}]\label{de12}
We say that $\Omega$ satisfies the exterior $C^{1}$ condition at $x_0\in \partial \Omega$
if there exist $r_0>0$ and a coordinate system $\{x_1,...,x_n \}$ (isometric to the original coordinate system) such that $x_0=0$ in this coordinate system and
\begin{equation}\label{deeq12}
B_{r_0} \cap \{(x',x_n)\big|x_n <-|x'|\omega(|x'|)\} \subset B_{r_0}\cap \Omega^c,
\end{equation}
where $\omega$ is a modulus of continuity. In addition, if $\omega$ is a Dini function, we say that  $\Omega$ satisfies the exterior $C^{1,\mathrm{Dini}}$ condition at $x_0$.
\end{definition}

\begin{remark}\label{r-1.2}
If $\Omega$ satisfies the exterior $C^{1, \mathrm{Dini}}$ condition at every $x_0\in \partial \Omega$ with the same $r_0$ and $\omega$, $\Omega$ is called a $\gamma$-convex domain (see \cite[Definition 1.1]{MR3029135}). If $\omega\equiv 0$, $\Omega$ is just a convex domain.
\end{remark}

\begin{definition}[\textbf{interior ${C^1}$ condition}]\label{de14}
We say that $\Omega$ satisfies the interior ${C^1}$ condition at ${x_0}\in \partial\Omega$ if there exist ${r_0}>0$
and a coordinate system $\left\{{{x_1},\ldots ,{x_n}}\right\}$ (isometric to the original coordinate system) such that ${x_0}=0$ in this coordinate system and
\begin{equation}\label{deeq14}
B_{r_0} \cap \{(x',x_n)\big|x_n >|x'|\omega(|x'|)\} \subset B_{r_0}\cap \Omega,
\end{equation}
where $\omega$ is a modulus of continuity.
\end{definition}

We also need the following definition for the pointwise behavior of a function.
\begin{definition}\label{de15}
Let $\Omega\subset {R^n}$ be a bounded set (may be not a domain) and $f$ be a function defined on $\Omega$.
We say that $f$ is differentiable at $x_0\in \Omega$ if there exist $r_0>0$, $l\in R^n$ and a modulus of continuity $\omega$ such that
\begin{equation*}
  |f(x)-f(x_0)-l\cdot (x-x_0)|\leq |x-x_0|\omega(|x-x_0|),~~\forall~x\in \Omega\cap B_{r_0}(x_0).
\end{equation*}
Then we denote $l$ by $\nabla f(x_0)$ and define
\begin{equation*}
  \|f\|_{C^{1}(x_0)}=|f(x_0)|+|\nabla f(x_0)|.
\end{equation*}
If $\omega$ is a Dini function, we say that $f$ is $C^{1,\mathrm{Dini}}$ at $x_0$ or $f\in C^{1,\mathrm{Dini}}(x_0)$.

Furthermore, $f$ is called $C^{-1,\mathrm{Dini}}$ at $x_0$ or $f\in C^{-1,\mathrm{Dini}}(x_0)$ if there exist $r_0>0$ and a Dini function $\omega$ such that
\begin{equation*}
\|f\|_{L^{n}(\Omega\cap B_{r}(x_0))}\leq \omega(r),~\forall ~0<r<r_0.
\end{equation*}

Finally, if $f\in C^{1,\mathrm{Dini}}(x_0)$ (or $f\in C^{-1,\mathrm{Dini}}(x_0)$) for any $x_0\in \Omega$ with the same $r_0$ and $\omega$, we say that $f\in C^{1,\mathrm{Dini}}(\Omega)$ (or $f\in C^{-1,\mathrm{Dini}}(\Omega)$).
\end{definition}

\begin{remark}\label{re1.3}
In this paper, we only use above definition when $\Omega$ is the closure of a bounded domain or $\Omega$ is the part of the boundary of a domain.
\end{remark}

\begin{remark}\label{re1.4}
Without loss of generality, we always assume that $r_0=1$ in \Crefrange{de12}{de15} throughout this paper.
\end{remark}

\begin{remark}\label{re1.2}
For a domain $\Omega$, we will use $\sigma_{\Omega}$ (resp., $\omega_{\Omega}$) to denote the corresponding modulus of continuity (resp., Dini function) in \Cref{de12} and \Cref{de14}. Similarly, for a function $f$, we will use $\sigma_{f}$ (resp., $\omega_{f}$) to denote the corresponding modulus of continuity (resp., Dini function) in \Cref{de15}.
\end{remark}

\subsection{Basic lemmas}\label{SP-2}

In this subsection, we introduce some basic lemmas which will be used in the proofs of our main results. The crucial one is the boundary pointwise $C^{1,\alpha}$ regularity at flat boundaries. The boundary $C^{1,\alpha}$ regularity was first proved by Krylov \cite{MR688919} for classical solutions and further simplified by Caffarelli (see \cite[Theorem 9.31]{MR1814364} and \cite[Theorem 4.28]{MR787227}).

For completeness, we give the detailed proof of the pointwise $C^{1,\alpha}$ regularity at a flat boundary by showing that the solution is controlled by two hyperplanes and is closer to one of them with aid of the Harnack inequality. We need the following simple result to prove the boundary pointwise $C^{1,\alpha}$ regularity.

\begin{lemma}\label{l.bas}
Let $u(e_n/2)\geq 1$ and $u\geq 0$ satisfy
\begin{equation*}
  \left\{\begin{aligned}
    &u\in S(\lambda,\Lambda,0)&& ~~\mbox{in}~~B_1^+;\\
    &u=0&& ~~\mbox{on}~~T_1.
  \end{aligned}\right.
\end{equation*}
Then
\begin{equation*}\label{e.bas}
u(x)\geq cx_n ~~\mbox{in}~~B^+_{1/2},
\end{equation*}
where $0<c<1$ is a universal constant.
\end{lemma}

\proof By the Harnack inequality \cite[Theorem 4.3]{MR1351007}, $u\geq c_0$ on $\partial B(e_n/4,1/8)$ for some positive universal constant $c_0$. Let
\begin{equation*}
v(x)=c_1\left(\left|x-e_n/4\right|^{-\alpha}-\left(1/4\right )^{-\alpha}\right),
\end{equation*}
where $c_1$ is small and $\alpha$ is large such that
\begin{equation*}
\left\{\begin{aligned}
  &\mathcal{M}^-(D^2v)\geq 0&&~~~~\mbox{in}~~~~B(e_n/4,1/4)\backslash \bar{B}(e_n/4,1/8);\\
  &v=0&&~~~~\mbox{on}~~~~\partial B(e_n/4,1/4);\\
  &v\leq c_0&&~~~~\mbox{on}~~~~\partial B(e_n/4,1/8).
\end{aligned}\right.
\end{equation*}

By the comparison principle and a direct calculation, we have
\begin{equation*}\label{e.bas-res1}
u\geq v\geq cx_n~~\mbox{on}~~\{x'=0,0\leq x_n\leq 1/8\}.
\end{equation*}
By translating $v$ to proper positions (with $v(x', 0)=0$ for some $x'\in T_{1/2}$) and similar arguments, we obtain
\begin{equation*}\label{e.bas-res2}
u\geq v\geq cx_n~~\mbox{in}~~\left\{|x'|<1/2,0<x_n<1/8\right\}.
\end{equation*}
Hence, by the Harnack inequality again,
\begin{equation*}\label{e.bas-res}
u\geq cx_n~~\mbox{in}~~B^+_{1/2}.
\end{equation*}
\qed~\\

Now, we can prove the boundary pointwise $C^{1,\alpha}$ regularity at a flat boundary based on above lemma.
\begin{lemma}\label{le21}
Let $u$ be a viscosity solution of
\begin{equation*}
\left\{\begin{aligned}
&u\in S(\lambda,\Lambda,0)&& ~~\mbox{in}~~B_1^+;\\
&u=0&& ~~\mbox{on}~~T_1.
\end{aligned}\right.
\end{equation*}

Then $u$ is $C^{1,{\alpha_1}}$ at $0$, i.e., there exists a constant $a$ such that
\begin{equation}\label{e.reg-res}
  |u(x)-ax_n|\leq C |x|^{1+{\alpha_1}}\|u\|_{L^{\infty }(B_1^+)}, ~~\forall ~x\in B_{1}^+
\end{equation}
and
\begin{equation*}
  |a|\leq C\|u\|_{L^{\infty }(B_1^+)},
\end{equation*}
where $0<\alpha_1<1$ and $C$ are universal constants.
\end{lemma}

\begin{remark}\label{re2.1}
From now on, $\alpha_1$ is fixed.
\end{remark}

\proof Without loss of generality, we assume that $\|u\|_{L^{\infty}(B_1^+)}=1$. To show \cref{e.reg-res}, we only need to prove the following: there exist a nonincreasing sequence $\{a_k\}$ ($k\geq 0$) and a nondecreasing sequence $\{b_k\}$ ($k\geq 0$) such that for all $k\geq 1$,
\begin{equation}\label{e.reg-M}
\begin{aligned}
  &b_kx_n\leq u\leq a_kx_n~~~~\mbox{in}~~~~B^+_{2^{-k}},\\
  &0\leq a_k-b_k\leq (1-\mu)(a_{k-1}-b_{k-1}),
\end{aligned}
\end{equation}
where $0<\mu<1/2$ is universal.

We prove the above by induction. Let $v(x)=2(1-|x+e_n|^{-n\Lambda/\lambda})$. Then $v$ satisfies
\begin{equation*}
\left\{\begin{aligned}
  &\mathcal{M}^+(D^2v)\leq 0&&~~~~\mbox{in}~~~~B^+_1;\\
  &v\geq 0 &&~~~~\mbox{in}~~~~B^+_1;\\
  &v\geq 1 &&~~~~\mbox{on}~~~~\partial B_1\backslash T_1;\\
  &v(0)=0.
\end{aligned}\right.
\end{equation*}
By the comparison principle and a direct calculation,
\begin{equation*}
-Cx_n\leq -v\leq u\leq v\leq Cx_n ~~\mbox{on}~~\{x'=0,0\leq x_n\leq 1/2\}.
\end{equation*}
By translating $v$ to proper positions (with $v(x', 0)=0$ for some $x'\in T_{1/2}$) and similar arguments,
\begin{equation*}
-Cx_n\leq u\leq Cx_n ~~\mbox{in}~~ B^+_{1/2}.
\end{equation*}
Thus, by taking $a_1=C, b_1=-C$ and $a_0=2C, b_0=-2C$, \cref{e.reg-M} holds for $k=1$ where $0<\mu<1/2$ is to be specified later.

Assume that \cref{e.reg-M} holds for $k$ and we need to prove it for $k+1$. Since \cref{e.reg-M} holds for $k$, there are two possible cases:
\begin{equation*}
\begin{aligned}
&\mbox{\textbf{Case 1}:}~~&&~~u(2^{-k-1}e_n)\geq \frac{a_k+b_k}{2}\cdot \frac{1}{2^{k+1}},\\
&\mbox{\textbf{Case 2}:}~~&&~~u(2^{-k-1}e_n)< \frac{a_k+b_k}{2}\cdot \frac{1}{2^{k+1}}.
\end{aligned}
\end{equation*}
Without loss of generality, we suppose that \textbf{Case 1} holds. Let $w=u-b_kx_n$ and then
\begin{equation*}
  \left\{\begin{aligned}
    &w\in S(\lambda,\Lambda,0)&& ~~\mbox{in}~~B_{2^{-k}}^+;\\
    &w\geq 0&& ~~\mbox{in}~~B_{2^{-k}}^+;\\
    &w(2^{-k-1}e_n)\geq \frac{a_k-b_k}{2}\cdot \frac{1}{2^{k+1}}.
  \end{aligned}\right.
\end{equation*}

From \Cref{l.bas} (the scaling version), for some universal constant $0<\mu<1$,
\begin{equation*}\label{e.reg-est}
w(x)\geq \mu(a_k-b_k)x_n~~\mbox{in}~~B_{2^{-k-1}}^+.
\end{equation*}
Thus,
\begin{equation}\label{e.reg-est-u}
u(x)\geq (b_k+\mu(a_k-b_k))x_n~~\mbox{in}~~B_{2^{-k-1}}^+.
\end{equation}
Let $a_{k+1}=a_k$ and $b_{k+1}=b_k+\mu(a_k-b_k)$. Then
\begin{equation}\label{e.reg-ak}
a_{k+1}-b_{k+1}=(1-\mu)(a_k-b_k).
\end{equation}
Hence, \cref{e.reg-M} holds for $k+1$. By induction, the proof is completed. \qed~\\

The next is a Hopf type lemma.
\begin{lemma}\label{le-1}
Let $u(e_n/2)\geq 1$ and $u\geq 0$ satisfy
\begin{equation*}
  u\in S^{\ast}(\lambda,\Lambda,f)~~\mbox{in}~~B_1^+.
\end{equation*}
Then there exists a universal constant $0<\delta_0<1$ such that if $\|f\|_{L^n(B_1^+)}\leq \delta_0$,
\begin{equation}\label{e-H-l}
u(x)\geq cx_n-C\|f\|_{L^n(B_1^+)} ~~\mbox{in}~~B_{1/2}^+,
\end{equation}
where $c>0$ and $C$ are universal constants.
\end{lemma}
\proof By the interior Harnack inequality, there exists a universal constant $0<\delta_0<1$ such that
\begin{equation*}
  \inf_{B(e_n/4,1/8)} u\geq 2\delta_0u(e_n/2)-\|f\|_{L^n(B_1^+)}\geq \delta_0.
\end{equation*}
Take
\begin{equation*}
v(x)=c_1\left(\left|x-e_n/4\right|^{-\alpha}-\left(1/4\right )^{-\alpha}\right).
\end{equation*}
As in the proof of \Cref{l.bas}, by choosing $c_1$ small and $\alpha$ large enough,
\begin{equation*}
\left\{\begin{aligned}
  &\mathcal{M}^-(D^2v)\geq 0&&~~~~\mbox{in}~~~~B(e_n/4,1/4)\backslash \bar{B}(e_n/4,1/8);\\
  &v=0&&~~~~\mbox{on}~~~~\partial B(e_n/4,1/4);\\
  &v\leq \delta_0&&~~~~\mbox{on}~~~~\partial B(e_n/4,1/8).
\end{aligned}\right.
\end{equation*}
Hence,
\begin{equation*}\label{e.bas-res1}
u\geq v\geq cx_n~~\mbox{on}~~\{x'=0,0\leq x_n\leq 1/8\}.
\end{equation*}
Let $w=u-v$ and then
\begin{equation*}
  \left\{\begin{aligned}
    \mathcal{M}^-(D^2w)&\leq f~~\mbox{in}~B(e_n/4,1/4)\backslash B(e_n/4,1/8);\\
    w&\geq 0~~\mbox{on}~\partial\big(B(e_n/4,1/4)\backslash B(e_n/4,1/8)\big).
\end{aligned}\right.
\end{equation*}
By the Alexandrov-Bakel'man-Pucci maximum principle \cite[Theorem 3.2]{MR1351007},
\begin{equation*}
  w\geq -C\|f\|_{L^n(B_1^+)}~~~~\mbox{in}~B(e_n/4,1/4)\backslash B(e_n/4,1/8).
\end{equation*}
Therefore,
\begin{equation*}
  u=v+w\geq cx_n-C\|f\|_{L^n(B_1^+)}~~~~\mbox{on}~\left\{x: |x'|=0, 0\leq x_n\leq 1/8\right\}.
\end{equation*}

By translating $v$ to proper positions (with $v(x', 0)=0$ for some $x'\in T_{1/2}$) and similar arguments, we have
\begin{equation*}
  u\geq cx_n-C\|f\|_{L^n(B_1^+)}~~~~\mbox{in}~\left\{x: |x'|<1/2, 0\leq x_n\leq 1/8\right\}.
\end{equation*}
Apply the Harnack inequality again and we derive
\begin{equation*}
 u\geq cx_n-C\|f\|_{L^n(B_1^+)}~~~~\mbox{in}~B_{1/2}^+.
\end{equation*}
\qed~\\

The following lemma is a simple consequence of the strong maximum principle and \Cref{l.bas}.
\begin{lemma}\label{le22}
Let $\Gamma \subset \partial {B_1^+} \backslash T_1$ and $u$ satisfy
\begin{equation*}
\left\{\begin{aligned}
&u\in S(\lambda,\Lambda,0) && ~~\mbox{in}~~B_1^+;\\
&u=x_n && ~~\mbox{on}~~\Gamma;\\
&u=0 && ~~\mbox{on}~~\partial {B_1^+}\backslash\Gamma.
\end{aligned}\right.
\end{equation*}
Then
\begin{equation}\label{e.2.1}
  u\ge cx_n~~\mbox{in}~~B^+_{1/2},
\end{equation}
where $0<c<1$ depends only on $n, \lambda$, $\Lambda$ and $\Gamma$.
\end{lemma}
\proof By the strong maximum principle and noting  $u\geq 0$, $u(e_n/2)\geq c_0>0$ where $c_0$ depends only on $n,\lambda,\Lambda$ and $\Gamma$. Then from \Cref{l.bas}, there exists $c>0$ such that \cref{e.2.1} holds.~\qed~\\

\section{Proof of \Cref{th11}}\label{Sth11}

The main result proved in this section is \Cref{th11}. Firstly, we prove a special case of \Cref{th11}, which shows the key idea clearly for the boundary differentiability.
\begin{theorem}\label{th21}
Let the assumptions of \Cref{th11} be satisfied. Suppose further that there exists a modulus of continuity $\sigma$ such that
\begin{equation}\label{e.ass-2}
\begin{aligned}
&\|u\|_{L^{\infty}(\Omega_1)}\leq c_0,~\|f\|_{L^n(\Omega_r)}\leq c_0\sigma(r), ~\forall ~0<r<1,\\
&-|x'|\sigma(|x'|)\leq x_n,~\forall~x\in (\partial \Omega)_1,\\
& |g(x)| \leq c_0|x|\sigma(|x|),~\forall~x\in (\partial \Omega)_1,\\
&\sigma(1)\leq c_0~~\mbox{and}~~
r^{\beta}\int_{r}^{1}\frac{\sigma(t)}{t^{1+\beta}}dt\leq c_0, ~\forall ~0<r<1,
\end{aligned}
\end{equation}
where $0<c_0<1/4$ and $0<\beta<1$ are constants (depending only on $n,\lambda,\Lambda$ and $K$) to be specified later.

Then there exist a universal constant $\bar{C}$ and a nonnegative sequence of $\{a_k\}$ $(k\ge -1)$ such that for all $k\geq 0$,
\begin{equation}\label{eq22}
\sup_{\Omega _{\eta^{k}}}(u-a_kx_n)\leq \eta ^{k}A_k,
\end{equation}
\begin{equation}\label{e.1.7}
\inf_{\Omega _{\eta^{k}}}(u+a_kx_n)\geq -\eta ^{k}A_k
\end{equation}
and
\begin{equation}\label{eq23}
{a_k} \le (1-\mu){a_{k-1}}+\bar C A_{k-1},
\end{equation}
where $0<\eta,\mu<1/4$ depend only on $n,\lambda,\Lambda$ and $K$,
\begin{equation}\label{eq24}
  A_{-1}=0,~A_0=c_0 ~~\mbox{and}~~A_k=\max(\sigma (\eta^{k-1}),~\eta^{\alpha_1/2}A_{k-1})~(k\ge 1).
\end{equation}
\end{theorem}
\begin{remark}\label{r3.1}
Note that $0<\alpha_1<1$ is the universal constant originated from \Cref{le21}.
\end{remark}
\proof We only give the proof of \Cref{eq22} and \Cref{eq23}. The proof of \cref{e.1.7} is similar and we omit it. We prove \Cref{eq22} and \Cref{eq23} by induction. For $k=0$, by setting $a_0=a_{-1}=0$, they hold clearly. Suppose that they hold for $k$. We need to prove that they hold for $k+1$.

Set $r=\eta ^{k}/2$, $\tilde{B}^{+}_{r}=B^{+}_{r}-r\sigma(r)e_n $, $\tilde{T}_r=T_r-r\sigma(r) e_n$ and
$\tilde{\Omega }_{r}=\Omega \cap \tilde{B}^{+}_{r}$.
Note that $\sigma(\eta) \leq \sigma(1)\leq c_0\leq 1/4$ and hence $\Omega _{r/2}\subset \tilde{\Omega }_{r}\subset \Omega_{2r}$. Let $v_1$ solve
\begin{equation*}
\left\{\begin{aligned}
 &\mathcal{M}^+({D^2}{v_1})=0 &&\mbox{in}~~\tilde{B}^{+}_{r}; \\
 &v_1=0 &&\mbox{on}~~\tilde{T}_{r};\\
 &v_1=\eta^kA_k &&\mbox{on}~~\partial \tilde{B}^{+}_{r}\backslash \tilde{T}_{r}.
\end{aligned}
\right.
\end{equation*}
By the boundary $C^{1,\alpha}$ estimate for $v_1$ (see \Cref{le21}) and the maximum principle, there exists $\bar a\ge 0$ such that (note that $\eta^{\alpha_1/2}A_k\leq A_{k+1}$)
\begin{equation*}
\begin{aligned}
 \|v_1-\bar{a}(x_n+r\sigma(r))\|_{L^{\infty }(\Omega_{2\eta r})}&\leq C\frac{(3\eta r)^{1+\alpha_1}}{r^{1 + \alpha_1 }}
 \|v\|_{L^\infty(\tilde{B}^{+}_{r})}\leq C\eta^{1+\alpha_1} \eta^kA_k
 \leq C\eta^{\alpha_1/2}\eta^{k+1}A_{k+1}
\end{aligned}
\end{equation*}
and
\begin{equation*}
\bar a \le \bar{C}A_{k}
  \end{equation*}
provided
\begin{equation}\label{e.1.4}
c_0\leq \eta,
\end{equation}
where $\bar{C}$ is universal. Hence,
\begin{equation}\label{eq28}
\begin{aligned}
\|v_1-\bar{a}x_n\|_{L^\infty(\Omega_{\eta^{k+1}})}=&\|v_1-\bar{a}x_n\|_{L^\infty(\Omega_{2\eta r})}\\
\leq & C\eta^{\alpha_1/2}\eta^{k+1}A_{k+1}+|\bar{a}r\sigma(r)|\\
\leq &\left(C\eta^{\alpha_1/2}+\frac{\bar{C}\sigma(\eta^k)}{\eta^{1+\alpha_1/2}}\right)
\cdot\eta^{k+1}A_{k+1}\\
\leq &\left(C\eta^{\alpha_1/2}+\frac{c_0C}{\eta^{1+\alpha_1/2}}\right)
\cdot\eta^{k+1}A_{k+1}.
\end{aligned}
\end{equation}

Let $v_2$ solve
\begin{equation*}
\left\{\begin{aligned}
 & \mathcal{M}^-(D^2{v_2})=0 &&\mbox{in}~~{{\tilde B_r^ +}}; \\
 & v_2=a_kx_n &&\mbox{on}~~{\partial {\tilde B_r^ +}\cap K};\\
 & v_2=0 &&\mbox{on}~~{\partial {\tilde B_r^ +}\backslash K}.
\end{aligned}
\right.
\end{equation*}
By \Cref{le22}, there exists $0<\mu<1/4$ (depending only on $n,\lambda$, $\Lambda$ and $K$) such that
\begin{equation}\label{eq26}
v_2\geq \mu a_k(x_n+r\sigma(r))\geq \mu a_kx_n~~~~\mbox{in}~  \Omega\cap B_{\eta^{k+1}}.
\end{equation}
In addition, it is easy to verify that
\begin{equation*}
v_2\leq a_k(x_n+r\sigma(r))~~~~\mbox{on}~\partial \tilde B_r^ +.
\end{equation*}
Then by the comparison principle,
\begin{equation*}
v_2\le a_k(x_n+r\sigma(r))~~~~\mbox{in}~ \tilde B_r^ +.
\end{equation*}

Let $w=u-a_kx_n-v_1+v_2$ and it follows that (note that $v_1,v_2\geq 0$)
\begin{equation*}
\left\{
    \begin{aligned}
&w\in \underline{S}(\lambda,\Lambda,-|f|) &&\mbox{in}~~ \Omega \cap \tilde{B}^{+}_{r}; \\
&w\leq (c_0+a_k)r\sigma(r) &&\mbox{on}~~\partial \Omega \cap \tilde{B}^{+}_{r};\\
&w\leq 0 &&\mbox{on}~~\partial \tilde{B}^{+}_{r}\cap \bar{\Omega}.
    \end{aligned}
    \right.
\end{equation*}
By the Alexandrov-Bakel'man-Pucci maximum principle, we have
\begin{equation*}
\begin{aligned}
\sup_{\Omega_{\eta^{k+1}}}w\leq \sup_{\tilde{\Omega}_{r}}w
& \leq (c_0+a_k)r \sigma(r)+Cr\|f\|_{L^n(\tilde{\Omega}_{r})}\leq (2c_0C+a_k)\eta^k \sigma(\eta^k).
\end{aligned}
\end{equation*}

Now, we try to estimate $a_k$. From the definition of $a_k$ (see \cref{eq23}),
\begin{equation}\label{e.1.6}
  \begin{aligned}
    a_k&\leq \bar{C}\sum_{i=0}^{k-1}(1-\mu)^{k-1-i}A_i\leq \bar{C}(1-\mu)^{k-1}\sum_{i=0}^{k-1}(1-\mu)^{-i}A_i.
  \end{aligned}
\end{equation}
From the definition of $A_i$ (see \cref{eq24}), for any $i\geq 1$,
\begin{equation*}
A_i\leq \sigma(\eta^{i-1})+\eta^{\alpha_1/2}A_{i-1}
\end{equation*}
and then
\begin{equation*}
  \begin{aligned}
\sum_{i=1}^{k-1}(1-\mu)^{-i}A_i\leq \sum_{i=1}^{k-1}(1-\mu)^{-i}\sigma(\eta^{i-1})
+\sum_{i=0}^{k-1}\eta^{\alpha_1/2}(1-\mu)^{-i-1}A_i.
  \end{aligned}
\end{equation*}
Hence, if
\begin{equation}\label{e.1.5}
\eta^{\alpha_1/2}\leq (1-\mu)/2,
\end{equation}
we have
\begin{equation*}
  \begin{aligned}
\sum_{i=1}^{k-1}(1-\mu)^{-i}A_i
&\leq 2\sum_{i=1}^{k-1}(1-\mu)^{-i}\sigma(\eta^{i-1})+c_0\\
&= c_0+\frac{2\sigma(1)}{1-\mu}+\frac{2}{(1-\mu)^2(1-\eta)}
\sum_{i=1}^{k-2}(1-\mu)^{-i+1}\frac{\sigma(\eta^{i})}{\eta^{i-1}}(\eta^{i-1}-\eta^i)\\
&\leq 4c_0+\frac{2}{(1-\mu)^2(1-\eta)}\int_{\eta^k}^{1}\frac{\sigma(t)}{t^{1+\beta}}dt,
  \end{aligned}
\end{equation*}
where $0<\beta<1$ satisfies $1-\mu=\eta^{\beta}$. By combining with \cref{e.1.6}, we have (note that $0<\eta,\mu<1/4$)
\begin{equation}\label{e.1.8}
  \begin{aligned}
a_k&\leq \frac{\bar{C}}{1-\mu}\eta^{k\beta}
    \left(4c_0+\frac{2}{(1-\mu)^2(1-\eta)}\int_{\eta^k}^{1}\frac{\sigma(t)}{t^{1+\beta}}dt\right)\\
&=\frac{4c_0\bar{C}}{1-\mu}\eta^{k\beta}+\frac{2\bar{C}\eta^{k\beta}}{(1-\mu)^3(1-\eta)}
\int_{\eta^k}^{1}\frac{\sigma(t)}{t^{1+\beta}}dt\\
&\leq \frac{4c_0\bar{C}}{1-\mu}
+\frac{2c_0\bar{C}}{(1-\mu)^3(1-\eta)}
\leq c_0C.
  \end{aligned}
\end{equation}
Therefore, for $w$,
\begin{equation}\label{eq29}
\sup_{\Omega_{\eta^{k+1}}}w
\leq \frac{c_0C}{\eta}\cdot\eta^{k+1}A_{k+1}.
\end{equation}

Take $\eta$ small enough such that \cref{e.1.5} holds and
\begin{equation*}
C\eta^{\alpha_1/2}\le \frac{1}{3}.
\end{equation*}
Next, choose $c_0$ small enough such that \cref{e.1.4} holds and
\begin{equation*}
\frac{c_0C}{\eta^{1+\alpha_1/2}}\leq \frac{1}{3}.
\end{equation*}
Let $a_{k+1}=(1-\mu)a_k+\bar a$. Then \Cref{eq23} holds for $k+1$. Recalling \Cref{eq28}, \Cref{eq26} and \Cref{eq29}, we have
\begin{equation*}
  \begin{aligned}
u-a_{k+1}x_n&=u-a_kx_n-v_1+v_2+v_1-\bar ax_n+\mu a_kx_n-v_2\\
&=w+v_1-\bar ax_n+\mu a_kx_n-v_2\\
&\le w+v_1-\bar ax_n\\
&\le \eta^{k+1}A_{k+1}~~\mbox{in} ~~\Omega_{\eta^{k+1}}.
  \end{aligned}
\end{equation*}
By induction, the proof is completed.\qed~\\

Now, we give the~\\
\noindent\textbf{Proof of \Cref{th11}.} In the following, we only need to show two things: One is that by a proper transformation, the assumptions \cref{e.ass-2} can be guaranteed. Another is that \crefrange{eq22}{eq23} imply the differentiability of $u$ at $0$.

For $\rho >0$, let $y=x/\rho$, $C=\|u\|_{L^{\infty}(\Omega_1)}+\|g\|_{C^{1}(0)}+1$ and
\begin{equation*}
  v(y)=\frac{c_0}{C}\left(u(x)-g(0)-\nabla g(0)\cdot x\right),
\end{equation*}
where $c_0$ is as in \cref{e.ass-2}. Then $v$ satisfies
\begin{equation*}
\left\{\begin{aligned}
&v\in S^*(\lambda,\Lambda,\tilde{f})&&~~\mbox{in}~~\tilde{\Omega}\cap B_1;\\
&v=\tilde{g}&& ~~\mbox{on}~~\partial \tilde{\Omega}\cap B_1,
\end{aligned}\right.
\end{equation*}
where
\begin{equation*}
\tilde{f}(y)=\frac{c_0\rho^{2}}{C}f(x),
  ~~\tilde{g}(y)=\frac{c_0}{C}\left(g(x)-g(0)-\nabla g(0)\cdot x\right)~~\mbox{and}
  ~~\tilde{\Omega}=\frac{\Omega}{\rho}.
\end{equation*}

First, it can be verified easily that
\begin{equation*}
  \begin{aligned}
&\|v\|_{L^{\infty}(\tilde{\Omega}_1)}\leq c_0,\\
&\|\tilde{f}\|_{L^{n}(\tilde{\Omega}_r)}=\frac{c_0\rho}{C}\|f\|_{L^{n}(\Omega_{\rho r})}
\leq c_0\rho\sigma_f(\rho r):=c_0\sigma_{\tilde{f}}(r),~\forall ~0<r<1,\\
&|\tilde{g}(y)|=|\frac{c_0}{C}(g(x)-g(0)-\nabla g(0)\cdot x)|\leq
c_0|x|\sigma_g(|x|)=c_0\rho |y|\sigma_g(\rho |y|):=c_0|y|\sigma_{\tilde{g}}(|y|),\\
  \end{aligned}
\end{equation*}
Next, for $x\in (\partial \Omega)_{1}$,
\begin{equation*}
-|x'|\sigma_{\Omega}(|x'|)\leq x_n,
\end{equation*}
where $\sigma_{\Omega}$ is the modulus of continuity corresponding to the exterior $C^{1}$ condition at $0$. Define $\sigma_{\tilde{\Omega}}(r)=\sigma_{\Omega}(\rho r)$ ($0<r<1$) and we have
\begin{equation*}
-|y'|\sigma_{\tilde{\Omega}}(|y'|)\leq y_n,~\forall~y\in (\partial \tilde\Omega)_1.
\end{equation*}
Let $\tilde\sigma(r)=\max(\sigma_{\tilde{f}}(r),\sigma_{\tilde{g}}(r),\sigma_{\tilde{\Omega}}(r))$ ($0<r<1$). Then
\begin{equation*}
  \tilde\sigma(1)\leq \sigma_f(\rho)+\sigma_g(\rho)+\sigma_{\Omega}(\rho)\rightarrow 0 ~\mbox{as}~\rho\rightarrow 0
\end{equation*}
and
\begin{equation*}
r^{\beta}\int_{r}^{1}\frac{\tilde\sigma(t)}{t^{1+\beta}}dt
\leq (\rho r)^{\beta}\int_{\rho r}^{1}
\frac{\sigma_f(t)+\sigma_g(t)+\sigma_{\Omega}(t)}{t^{1+\beta}}dt \rightarrow 0 ~\mbox{as}~\rho\rightarrow 0.
\end{equation*}
Therefore, by choosing $\rho$ small enough (depending only on $n,\lambda,\Lambda,K$, $\sigma_{f}$, $\sigma_{g}$ and $\sigma_{\Omega}$), the assumptions \cref{e.ass-2} for $v$ can be guaranteed. Hence, we can apply \Cref{th21} to $v$.

Finally, we show that the differentiability of $v$ (and hence $u$) at $0$ can be inferred by \crefrange{eq22}{eq23}. From \cref{eq23}, we know that $a_k\rightarrow 0$ as $k\rightarrow \infty$ (see also \cref{e.1.8}). For any $y\in \tilde\Omega_1$, there exists $k\geq 0$ such that $\eta^{k+1}\leq |y|<\eta^{k}$ and then define $\sigma_v(|y|)=(a_k+A_k)/\eta$. Thus, $\sigma_v(r)\rightarrow 0$ as $r\rightarrow 0$. From \Cref{eq22} and \cref{e.1.7}, we obtain
\begin{equation*}
  |v(y)|\leq \eta^{k}(a_k+A_k)\leq |y|\sigma_v(|y|).
\end{equation*}
Therefore, $v$ is differentiable at 0 with $\nabla v(0)=0$. By rescaling back to $u$, we obtain the conclusion of \Cref{th11}.\qed~\\

\section{Proof of \Cref{th12}}\label{Sth12}

In this section, we prove \Cref{th12}. First, we show that the graph of the solution can be controlled by two hyperplanes.
\begin{lemma}\label{le33}
Let the assumptions of \Cref{th12} be satisfied. Suppose further that there exists a Dini function $\omega$ such that
\begin{equation}\label{e.ass-0}
\begin{aligned}
&\|u\|_{L^{\infty}(\Omega_1)}\leq c_0,~\|f\|_{L^n(\Omega_r)}\leq c_0\omega(r), ~\forall ~0<r<1,\\
&-|x'|\omega(|x'|)\leq x_n,~\forall~x\in (\partial \Omega)_1,\\
& |g(x)| \leq c_0|x|\omega(|x|),~\forall~x\in (\partial \Omega)_1,\\
&\omega(1)\leq c_0~~\mbox{and}~~\int_{0}^{1}\frac{\omega(\rho)}{\rho}d\rho\leq c_0,
\end{aligned}
\end{equation}
where $0<c_0<1/4$ is a universal constant to be specified later.

Then there exist sequences of constants $\{\tilde{a}_k\},\{\tilde{b}_k\}$ and $\{\tilde{A}_k\}$ $(k\geq -1)$ with $\tilde{a}_k\geq 0\geq \tilde{b}_k$ and $\tilde{A}_k\geq 0$ such that for any $k\geq 0$,
\begin{equation}\label{e.eq31}
 \tilde{b}_kx_n -\eta^k\tilde{A}_k\leq u\leq \tilde{a}_kx_n+\eta^k\tilde{A}_k~~\mbox{in}~~\Omega_{\eta^k}
\end{equation}
and
\begin{equation}\label{e.s1-0}
\begin{aligned}
\tilde{a}_{k}\leq \tilde{a}_{k-1}+\bar{C}\tilde{A}_{k-1},~
\tilde{b}_{k}\geq\tilde{b}_{k-1}-\bar{C}\tilde{A}_{k-1},
  \end{aligned}
\end{equation}
where $0<\eta<1/4$ and $\bar{C}$ are universal constants,
\begin{equation*}
\tilde{A}_{-1}=0,~ \tilde{A}_0=c_0,~
\tilde{A}_{k}=\max(\omega(\eta^{k-1}), \eta^{\alpha_1/2}\tilde{A}_{k-1}) (k\geq 1).
\end{equation*}
\end{lemma}
\proof We prove \cref{e.eq31} and \cref{e.s1-0} by induction. For $k=0$, by setting $\tilde{a}_{-1}=\tilde{a}_0=\tilde{b}_{-1}=\tilde{b}_0=0$, the conclusion holds clearly. Suppose that the conclusion holds for $k$. We need to prove that the conclusion holds for $k+1$.

Set $r=\eta ^{k}/2$, $\tilde{B}^{+}_{r}=B^{+}_{r}-r\omega(r)e_n $, $\tilde{T}_r=T_r-r\omega(r) e_n$ and $\tilde{\Omega }_{r}=\Omega \cap \tilde{B}^{+}_{r}$. Note that $\omega(\eta) \leq \omega(1)\leq c_0\leq 1/4$ and then $\Omega _{r/2}\subset \tilde{\Omega }_{r}\subset \Omega_{2r}$. Let $v$ solve
\begin{equation*}
\left\{\begin{aligned}
 &\mathcal{M}^+(D^2v)=0 &&\mbox{in}~~\tilde{B}^{+}_{r}; \\
 &v=0 &&\mbox{on}~~\tilde{T}_{r};\\
 &v=\eta ^{k} \tilde{A}_k &&\mbox{on}~~\partial \tilde{B}^{+}_{r}\backslash \tilde{T}_{r}
\end{aligned}
\right.
\end{equation*}
and $w=u- \tilde{a}_kx_n-v$. Then $w$ satisfies (note that $v\geq 0$ in $\tilde{B}^{+}_{r}$)
\begin{equation*}
\left\{
\begin{aligned}
&w\in \underline{S}(\lambda,\Lambda , -|f|) &&\mbox{in}~~ \Omega \cap \tilde{B}^{+}_{r}; \\
&w\leq g- \tilde{a}_kx_n &&\mbox{on}~~\partial \Omega \cap \tilde{B}^{+}_{r};\\
&w\leq 0 &&\mbox{on}~~\partial \tilde{B}^{+}_{r}\cap \bar{\Omega}.
    \end{aligned}
    \right.
\end{equation*}

In the following arguments, we estimate $v$ and $w$ respectively. For $v$, similar to the estimate in \Cref{th21}, there exists $\bar{a}\geq 0$  such that
\begin{equation}\label{e.19}
\bar{a}\leq \bar{C} \tilde{A}_k
\end{equation}
and
\begin{equation}\label{e.20}
\begin{aligned}
\|v-\bar{a}x_n\|_{L^{\infty }(\Omega _{\eta^{k+1}})}
\leq \left( C\eta ^{\alpha_1/2}+\frac{c_0 C}{\eta^{1+\alpha_1/2}}\right)\cdot
\eta ^{k+1} \tilde{A}_{k+1}
\end{aligned}
\end{equation}
provided
\begin{equation}\label{e.1.4-2}
c_0\leq \eta.
\end{equation}

From \cref{e.s1-0},
\begin{equation*}
 -\bar{C}\sum_{i=0}^{\infty} \tilde{A}_i \leq \tilde{b}_k\leq \tilde{a}_k\leq \bar{C}\sum_{i=0}^{\infty} \tilde{A}_i,
 ~\forall ~k\geq 0.
\end{equation*}
Moreover, it is easy to show that
\begin{equation*}
  \sum_{i=0}^{\infty} \tilde{A}_i\leq \sum_{i=1}^{\infty}\omega(\eta^i)+\eta^{\alpha_1/2}\sum_{i=0}^{\infty} \tilde{A}_i+2c_0,
\end{equation*}
which implies
\begin{equation}\label{e.1.1}
\begin{aligned}
\sum_{i=0}^{\infty} \tilde{A}_i
&\leq \frac{1}{1-\eta^{\alpha_1/2}}\sum_{i=1}^{\infty}\omega(\eta^i)
+\frac{2c_0}{1-\eta^{\alpha_1/2}}\\
&=\frac{1}{\left(1-\eta^{\alpha_1/2}\right)\left(1-\eta\right)}\sum_{i=1}^{\infty}
\frac{\omega(\eta^i)\left(\eta^{i-1}-\eta^i\right)}{\eta^{i-1}}+\frac{2c_0}{1-\eta^{\alpha_1/2}}\\
&\leq \frac{1}{\left(1-\eta^{\alpha_1/2}\right)\left(1-\eta\right)}\int_{0}^{1}
\frac{\omega(r)dr}{r}+\frac{2c_0}{1-\eta^{\alpha_1/2}}\\
&\leq \frac{c_0}{\left(1-\eta^{\alpha_1/2}\right)\left(1-\eta\right)}
+\frac{2c_0}{1-\eta^{\alpha_1/2}}\leq 6c_0,
\end{aligned}
\end{equation}
provided
\begin{equation}\label{e-w}
\left(1-\eta^{\alpha_1/2}\right)\left(1-\eta\right)\geq 1/2.
\end{equation}

For $w$, by the Alexandrov-Bakel'man-Pucci maximum principle, we have
\begin{equation}\label{e1.22}
  \begin{aligned}
\sup_{\Omega_{\eta^{k+1}}}w\leq \sup_{\tilde{\Omega} _{r}}w
& \leq\|g\|_{L^{\infty }(\partial \Omega \cap \tilde{B}^{+}_{r})}+\tilde{a}_kr\omega(r)
+Cr\|f\|_{L^n(\tilde{\Omega}_{r})}\\
&\leq c_0\eta^k\omega(\eta^k) +6c_0\bar{C}\eta^k\omega(\eta^k)
+Cc_0\eta^k\omega(\eta^k)\\
&\leq \frac{c_0C}{\eta}\cdot
\eta^{k+1}\tilde{A}_{k+1}.
  \end{aligned}
\end{equation}
Take $\eta $ small enough such that \cref{e-w} holds and
\begin{equation*}
  C\eta ^{\alpha_1/2}\leq \frac{1}{3}.
\end{equation*}
Take $c_0$ small enough such that \cref{e.1.4-2} holds and
\begin{equation*}
\frac{c_0C}{\eta^{1+\alpha_1/2}}\leq\frac{1}{3}.
\end{equation*}
Let $\tilde{a}_{k+1}=\tilde{a}_k+\bar{a}$ and \cref{e.s1-0} holds for $k+1$ clearly. By combining with \cref{e.20} and \cref{e1.22}, we have
\begin{equation*}
\begin{aligned}
u-\tilde{a}_{k+1}x_n&=u-\tilde{a}_kx_n-v+v-\bar{a}x_n=w+v-\bar{a}x_n\leq \eta ^{k+1}\tilde{A}_{k+1} ~~\mbox{in}~~\Omega_{\eta^{k+1}}.
\end{aligned}
\end{equation*}

Similarly, we can prove that there exists a nonnegative constant $\bar{b}$ with $\bar{b}\leq \bar{C}\tilde{A}_k$ such that for $\tilde{b}_{k+1}=\tilde{b}_{k}-\bar{b}$,
\begin{equation*}
\begin{aligned}
u-\tilde{b}_{k+1}x_n\geq -\eta ^{k+1}\tilde{A}_{k+1} ~~\mbox{in}~~\Omega_{\eta^{k+1}}.
\end{aligned}
\end{equation*}
Therefore, the proof is completed by induction.\qed~\\

\begin{remark}\label{re1.1}
From \cref{e.s1-0} and \cref{e.1.1},
\begin{equation*}
-1\leq \tilde{b}_k\leq \tilde{a}_k\leq 1,~\forall ~k\geq 0.
\end{equation*}
Then one can show that $u$ is Lipschitz continuous at $0$. In fact, this is one of the main results in \cite{MR1} and the proof is also adopted there. For completeness, we give the detailed proof here.
\end{remark}

Similarly to the last section, we prove a special case of \Cref{th12} first.
\begin{theorem}\label{th33}
Let the assumptions of \Cref{th12} be satisfied. Suppose further that there exist a Dini function $\omega$ and a modulus of continuity $\sigma$ with $\sigma\geq \omega$ such that
\begin{equation}\label{e.ass}
\begin{aligned}
&\|u\|_{L^{\infty}(\Omega_1)}\leq c_0,~\|f\|_{L^n(\Omega_r)}\leq c_0\omega(r), ~\forall ~0<r<1,\\
&-|x'|\omega(|x'|)\leq x_n\leq |x'|\sigma(|x'|),~\forall~x\in (\partial \Omega)_1,\\
& |g(x)| \leq c_0|x|\omega(|x|),~\forall~x\in (\partial \Omega)_1,\\
&\omega(1)\leq\sigma(1)\leq c_0~~\mbox{and}~~
\int_{0}^{1}\frac{\omega(\rho)}{\rho}d\rho\leq c_0,
\end{aligned}
\end{equation}
where $0<c_0<1/4$ is the universal constant as in \Cref{le33}.

Then there exist sequences of constants $\{a_k\},\{b_k\},\{A_k\}$ and $\{B_k\}$ $(k\geq 0)$ with $a_k\geq b_k$ and $A_k,B_k\geq 0$ such that for any $k\geq 0$,
\begin{equation}\label{eq31}
 b_kx_n -\eta^kB_k\leq u\leq a_kx_n+\eta^kA_k,~~\mbox{in}~~\Omega_{\eta^k}
\end{equation}
where $0<\eta<1/4$ is the universal constant as in \Cref{le33}; $a_0=b_0=0$ and $A_0=B_0=c_0$; for any $k\geq 1$,
\begin{equation}\label{e.s1-2}
\left\{\begin{aligned}
A_{k}&=\max(\sigma(\eta^{k-1}),\eta^{\alpha_1/2}A_{k-1}),\\
B_{k}&=2\sigma(\eta^{k-1})+\eta^{\alpha_1/2}B_{k-1}, \\
a_{k}&\leq \min(a_{k-1}+\bar{C}A_{k-1},\tilde{a}_{k}),\\
b_{k}&\geq \max(b_{k-1}+\mu (a_{k-1}-b_{k-1})-\bar{C}B_{k-1},\tilde{b}_{k})
  \end{aligned}
  \right.
\end{equation}
or
\begin{equation}\label{e.s1-1}
\left\{\begin{aligned}
A_{k}&=2\sigma(\eta^{k-1})+\eta^{\alpha_1/2}A_{k-1},\\
B_{k}&=\max(\sigma(\eta^{k-1}),\eta^{\alpha_1/2}B_{k-1}),\\
a_{k}&\leq \min(a_{k-1}-\mu (a_{k-1}-b_{k-1})+\bar{C}A_{k-1},\tilde{a}_{k}),\\
b_{k}&\geq \max(b_{k-1}-\bar{C}B_{k-1},\tilde{b}_{k}).
  \end{aligned}
  \right.
\end{equation}
Here, $0<\mu<1/8$ and $\bar{C}$ are universal constants, and the constants $\tilde{a}_k$ and $\tilde{b}_k$ are as in \cref{e.eq31}.
\end{theorem}

\proof We prove the above by induction. For $k=0$, the conclusion holds clearly. Suppose that the conclusion holds for $k$. We need to prove that the conclusion holds for $k+1$.

By similar arguments to the proof of \Cref{le33}, there exist nonnegative constants $\bar{a}$ and $\bar{b}$ with $\bar{a}\leq \bar{C}A_k$ and $\bar{b}\leq \bar{C}B_k$ such that (note that $\eta<1/4$)
\begin{equation}\label{e.1.2}
\begin{aligned}
u-\bar{a}_{k+1}x_n\leq \eta ^{k+1}\bar{A}_{k+1} ~~\mbox{in}~~\Omega_{3\eta^{k+1}}
\end{aligned}
\end{equation}
and
\begin{equation}\label{e.1.2-2}
\begin{aligned}
u-\bar{b}_{k+1}x_n\geq -\eta ^{k+1}\bar{B}_{k+1} ~~\mbox{in}~~\Omega_{3\eta^{k+1}},
\end{aligned}
\end{equation}
where
\begin{equation*}
  \begin{aligned}
&\bar{a}_{k+1}=a_{k}+\bar{a},&&~\bar{b}_{k+1}=b_{k}-\bar{b},~\\
&\bar{A}_{k+1}=\max(\sigma(\eta^k), \eta^{\alpha_1/2}A_k),~&&
\bar{B}_{k+1}=\max(\sigma(\eta^k), \eta^{\alpha_1/2}B_k).
  \end{aligned}
\end{equation*}

In the following, we prove the conclusion for $k+1$ according to two cases.

\textbf{Case 1:} $u(re_n)\geq \left(a_{k}+b_{k}\right)r/2$ where $r=\eta^{k+1}$. Let
\begin{equation*}
v(x)=u(x)-\bar{b}_{k+1}x_n+r\bar{B}_{k+1}.
\end{equation*}
Then $v\geq 0$ in $\Omega_{2r}$ and $v(re_n)\geq \left(a_k-b_k\right)r/2+r\bar{B}_{k+1}$.

Next, we apply \Cref{le-1}. Indeed, note that $ B_{2r}^++2r\sigma(2r)e_n \subset \Omega_{2r}$. For $x\in B_{2r}^++2r\sigma(2r)e_n$, set $y=(x-2r\sigma(2r)e_n)/(2r)$ and
\begin{equation*}
  w(y)=\frac{v(x)}{(a_k-b_k)r/2+r\bar{B}_{k+1}}.
\end{equation*}
Then $w\geq 0$ satisfies $w(y^*)\geq 1$ for some $y^*_n\geq 1/4$ and
\begin{equation*}
    \mathcal{M}^-(D^2w)\leq |\tilde{f}|~~\mbox{in}~~B_1^+,
\end{equation*}
where
\begin{equation*}
  \tilde{f}(y)=\frac{rf(x)}{(a_k-b_k)/2+\bar{B}_{k+1}}.
\end{equation*}
Note that
\begin{equation*}
  \|\tilde{f}\|_{L^n(B_1^+)}
  =\frac{\|f\|_{L^n(B_{2r}^++2r\sigma(2r)e_n)}}{(a_k-b_k)/2+\bar{B}_{k+1}}
  \leq \frac{c_0\omega(\eta^k)}{(a_k-b_k)/2+\bar{B}_{k+1}}\leq c_0.
\end{equation*}
Take $c_0\leq \delta_0$ (as in \Cref{le-1}). Then apply \Cref{le-1} to $w$, there exists a universal constant $0<\mu<1/8$ such that
\begin{equation*}
w\geq 2\mu y_n-C\|\tilde{f}\|_{L^n(B_1^+)} ~~\mbox{in}~~B_{1/2}^+.
\end{equation*}
By rescaling back to $v$, we have
\begin{equation*}
v\geq \mu\left(a_k-b_k\right)(x_n-2r\sigma(2r))-c_0Cr\omega(\eta^k) ~~\mbox{in}~~B_{r}^++2r\sigma(2r)e_n.
\end{equation*}
Note that $v\geq 0$ and thus the above holds in $\Omega_{r}$, i.e.,
\begin{equation*}
v\geq \mu\left(a_k-b_k\right)(x_n-2r\sigma(2r))-c_0Cr\omega(\eta^k) ~~\mbox{in}~~\Omega_{r}.
\end{equation*}
Hence,
\begin{equation}\label{e.1.3}
u\geq \left(\bar{b}_{k+1}+\mu\left(a_k-b_k\right)\right)x_n
-r\bar{B}_{k+1}-2\mu r(a_k-b_k)\sigma(\eta^k)-c_0Cr\omega(\eta^k)  ~~\mbox{in}~~\Omega_{r}.
\end{equation}
From \Cref{le33} (see also \Cref{re1.1}), $-1\leq\tilde{b}_k\leq\tilde{a}_k\leq 1$. Thus,
\begin{equation*}
  0\leq a_k-b_k\leq \tilde{a}_k-\tilde{b}_k\leq 2.
\end{equation*}
Let
\begin{equation*}
B_{k+1}=2\sigma(\eta^k)+\eta^{\alpha_1/2}B_k.
\end{equation*}
Then by taking $c_0$ small again (note that $\mu<1/8$),
\begin{equation*}
B_{k+1}\geq\bar{B}_{k+1}+2\mu(a_k-b_k)\sigma(\eta^k)+c_0C\omega(\eta^k).
\end{equation*}
Hence, by noting \cref{e.1.2} and \cref{e.1.3}, we have
\begin{equation}\label{e.2.2}
 \left(\bar{b}_{k+1}+\mu(a_k-b_k)\right)x_n -\eta^{k+1}B_{k+1}\leq u
 \leq \bar{a}_{k+1}x_n+\eta^{k+1}\bar{A}_k,~~\mbox{in}~~\Omega_{\eta^{k+1}}.
\end{equation}
Take
\begin{equation*}
  \begin{aligned}
&a_{k+1}=\min(\bar{a}_{k+1},\tilde{a}_{k+1})=\min(a_k+\bar{a},\tilde{a}_{k+1}),~\\
&b_{k+1}=\max(\bar{b}_{k+1}+\mu(a_k-b_k),\tilde{b}_{k+1})
=\max(b_{k}-\bar{b}+\mu(a_k-b_k),\tilde{b}_{k+1}),~\\
&A_{k+1}=\bar{A}_{k+1}=\max(\sigma(\eta^k), \eta^{\alpha_1/2}A_k).
  \end{aligned}
\end{equation*}
Thus, \cref{e.s1-2} holds for $k+1$. Note that $\tilde{a}_{k+1}\geq \tilde{a}_{k}\geq a_k$ and $\tilde{b}_{k+1}\leq \tilde{b}_{k}\leq b_k$. Then it can be checked easily that $a_{k+1}\geq b_{k+1}$. By combining with \cref{e.eq31}, $\tilde{A}_{k+1}\leq A_{k+1}$ and $\tilde{B}_{k+1}\leq B_{k+1}$, \cref{e.2.2} reduces to \cref{eq31} for $k+1$.

\textbf{Case 2:} $u(re_n)\leq \left(a_{k}+b_{k}\right)r/2$. Let $v(x)=\bar{a}_{k+1}x_n+r\bar{A}_{k+1}-u(x)$. Then by a similar argument (we omit it), \cref{eq31} and \cref{e.s1-1} hold for $k+1$. Therefore, the proof is completed by induction.\qed~\\

Now, we give the~\\
\noindent\textbf{Proof of \Cref{th12}.} As before, in the following, we only need to show two things: one is that by a proper transformation, the assumptions \cref{e.ass} can be guaranteed. Another is that \crefrange{eq31}{e.s1-1} imply the differentiability of $u$ at $0$.

For  $\rho>0$, let $y=x/\rho$, $C=\|u\|_{L^{\infty}(\Omega\cap B_1)}+\|g\|_{C^{1}(0)}+1$ and
\begin{equation*}
  v(y)=\frac{c_0}{C}\left(u(x)-g(0)-\nabla g(0)\cdot x\right).
\end{equation*}
Then $v$ satisfies
\begin{equation*}
\left\{\begin{aligned}
&v\in S^*(\lambda,\Lambda,\tilde{f})&&~~\mbox{in}~~\tilde{\Omega}\cap B_1;\\
&v=\tilde{g}&& ~~\mbox{on}~~\partial \tilde{\Omega}\cap B_1,
\end{aligned}\right.
\end{equation*}
where
\begin{equation*}
\tilde{f}(y)=\frac{c_0\rho^{2}}{C}f(x),
  ~~\tilde{g}(y)=\frac{c_0}{C}\left(g(x)-g(0)-\nabla g(0)\cdot x\right)~~\mbox{and}
  ~~\tilde{\Omega}=\frac{\Omega}{\rho}.
\end{equation*}

First, it can be verified easily that
\begin{equation*}
  \begin{aligned}
&\|v\|_{L^{\infty}(\tilde{\Omega}_1)}\leq c_0,\\
&\|\tilde{f}\|_{L^{n}(\tilde{\Omega}_r)}=\frac{c_0\rho}{C}\|f\|_{L^{n}(\Omega_{\rho r})}
\leq c_0\rho\omega_f(\rho r):=c_0\omega_{\tilde{f}}(r),\\
&|\tilde{g}(y)|=|\frac{c_0}{C}(g(x)-g(0)-\nabla g(0)\cdot x)|\leq
c_0|x|\omega_g(|x|)=c_0\rho|y|\omega_g(\rho |y|):=c_0|y|\omega_{\tilde{g}}(|y|).\\
  \end{aligned}
\end{equation*}
Next, for $x\in (\partial \Omega)_{1}$,
\begin{equation*}
-|x'|\omega_{\Omega}(|x'|)\leq x_n\leq |x'|\sigma_{\Omega}(|x'|),
\end{equation*}
where $\omega_{\Omega}$ is the Dini function corresponding to the exterior $C^{1,\mathrm{Dini}}$ condition at $0$ and $\sigma_{\Omega}$ is the modulus of continuity corresponding to the interior $C^{1}$ condition at $0$. Define $\sigma_{\tilde{\Omega}}(r)=\sigma_{\Omega}(\rho r)$ and $\omega_{\tilde{\Omega}}(r)=\omega_{\Omega}(\rho r)$ ($0<r<1$) and we have
\begin{equation*}
-|y'|\omega_{\tilde\Omega}(|y'|)\leq y_n\leq |y'|\sigma_{\tilde\Omega}(|y'|),~\forall~y\in (\partial \tilde\Omega)_1.
\end{equation*}
Let $\tilde\omega(r)=\max(\omega_{\tilde{f}}(r),\omega_{\tilde{g}}(r),\omega_{\tilde{\Omega}}(r))$ and $\tilde{\sigma}(r)=\max(\tilde{\omega}(r),\sigma_{\tilde\Omega}(r))$ ($0<r<1$) and then
\begin{equation*}
\tilde\sigma(1)
\leq \omega_f(\rho)+\omega_g(\rho)+\omega_{\Omega}(\rho)+\sigma_{\Omega}(\rho)
\rightarrow 0 ~\mbox{as}~\rho\rightarrow 0
\end{equation*}
and
\begin{equation*}
\int_{0}^{1}\frac{\tilde\omega(t)}{t}dt
\leq \int_{0}^{\rho}\frac{ \omega_f(t)+\omega_g(t)+\omega_{\Omega}(t)}{t}dt
\rightarrow 0 ~\mbox{as}~\rho\rightarrow 0.
\end{equation*}
Therefore, by choosing $\rho$ small enough (depending only on $n,\lambda,\Lambda$, $\omega_f$, $\omega_g$, $\omega_{\Omega}$ and $\sigma_{\Omega}$), the assumptions \cref{e.ass} for $v$ can be guaranteed. Hence, we can apply \Cref{th33} to $v$.

Next, we show that \crefrange{eq31}{e.s1-1} imply the differentiability of $v$ (and hence $u$) at $0$. Indeed, from \Cref{le33}, $\tilde{a}_k$ and $\tilde{b}_k$ are bounded sequences. By \cref{e.s1-2} and \cref{e.s1-1}, $a_k$ and $b_k$ are also bounded. Moreover, for any $k\geq 0$,
\begin{equation*}
  0\leq a_{k}-b_{k}\leq (1-\mu)(a_{k-1}-b_{k-1})+\tilde{C}(A_{k-1}+B_{k-1}).
\end{equation*}
Note that $A_k,B_k\rightarrow 0$ as $k\rightarrow \infty$. Hence, $a_{k}-b_{k}\rightarrow 0$ as $k\rightarrow \infty$.

In the following, we prove that the sequence $\{a_k\}$ converges, i.e., $\{a_k\}$ has at most one accumulation point. Suppose not and let $a$ and $\tilde{a}$ be two accumulation points of $\{a_k\}$. We assume that $a,\tilde{a}>0$ and $\varepsilon:=a-\tilde{a}>0$ without loss of generality. Since $A_k,B_k,a_k-b_k\rightarrow 0$ and $\omega$ is a Dini function, there exists $k_0\geq 0$ such that for any $k\geq k_0$,
\begin{equation*}
A_k,~B_k,~a_k-b_k\leq \varepsilon/8 ~~\mbox{and}~~\int_{0}^{\eta^{k}}\frac{\omega(r)dr}{r}\leq \frac{\varepsilon}{32\bar{C}}.
\end{equation*}
In addition, there exist $k_2\geq k_1\geq k_0$ such that
\begin{equation}\label{e.2.3}
a_{k_1}\leq \tilde{a}+\varepsilon/32 ~~\mbox{and}~~ a_{k_2}\geq a-\varepsilon/32.
\end{equation}

Since
\begin{equation*}
  u\leq a_{k_1}x_n+\eta^{k_1}A_{k_1}~~\mbox{in}~~\Omega_{\eta^{k_1}},
\end{equation*}
by a similar argument to the proof of \Cref{le33}, there exist a sequence of $\left\{\hat{a}_k\right\}_{k\geq k_1}$ such that for any $k\geq k_1$,
\begin{equation*}
  u\leq \hat{a}_kx_n+\eta^{k}\hat{A}_{k}~~\mbox{in}~~\Omega_{\eta^{k}}
\end{equation*}
and
\begin{equation*}
  \hat{a}_{k+1}\leq \hat{a}_k+\bar{C}\hat{A}_{k},
\end{equation*}
where
\begin{equation*}
  \hat{A}_{k_1}=A_{k_1}~~\mbox{and}~~
  \hat{A}_{k+1}=\max(\omega(\eta^k), \eta^{\alpha_1/2}\hat{A}_k) (k>k_1).
\end{equation*}
Thus, as before, for any $k\geq k_1$ (note \cref{e-w}),
\begin{equation*}
\begin{aligned}
\hat{a}_k\leq a_{k_1}+\bar{C}\sum_{i=k_1}^{\infty}\hat{A}_{i}
&\leq a_{k_1}+\frac{\bar{C}}{1-\eta^{\alpha_1/2}}
\sum_{i=k_1}^{\infty}\omega(\eta^i)+\frac{\bar{C}}{1-\eta^{\alpha_1/2}}A_{k_1}\\
&\leq a_{k_1}+\frac{\bar{C}}{\left(1-\eta^{\alpha_1/2}\right)\left(1-\eta\right)}
\int_{0}^{\eta^{k_1}}\frac{\omega(r)dr}{r}+\frac{\bar{C}}{1-\eta^{\alpha_1/2}}A_{k_1}\\
&\leq a_{k_1}+\frac{\varepsilon}{4}.
\end{aligned}
\end{equation*}
Then for $x\in \Omega_{\eta^{k_2}}$ with $x_n>0$,
\begin{equation*}
  \begin{aligned}
    &\left(a_{k_2}-\frac{\varepsilon}{8}\right)x_n-\eta^{k_2}B_{k_2}\leq b_{k_2}x_n-\eta^{k_2}B_{k_2}\leq u(x)\\
    &\leq \hat{a}_{k_1}x_n+\eta^{k_2}\hat{A}_{k_2}
    \leq \left(a_{k_1}+\frac{\varepsilon}{4}\right)x_n+\eta^{k_2}A_{k_2}.
  \end{aligned}
\end{equation*}
Choose $x=\eta^{k_2}e_n/2$ and we have
\begin{equation*}
  \begin{aligned}
a_{k_2}-a_{k_1}\leq \frac{7\varepsilon}{8},
  \end{aligned}
\end{equation*}
which contradicts with \cref{e.2.3}.

Therefore, there exists $\bar{a}$ such that $a_k,b_k\rightarrow \bar{a}$. Then \cref{eq31} implies that $v$ is differentiable at $0$. Indeed, for any $y\in \tilde\Omega_1$, there exists $k\geq 0$ such that $\eta^{k+1}\leq |y|<\eta^{k}$. By applying \cref{eq31} to $v$ in $\tilde\Omega_{\eta^k}$ and noting $a_k,b_k\rightarrow \bar{a}$ and $A_k,B_k\rightarrow 0$, we have
\begin{equation*}
  \begin{aligned}
    |v(y)-\bar{a}y_n|&\leq \max(|a_k-\bar{a}|,|b_k-\bar{a}|)|y|+\bar{C}(A_k+B_k)\eta^{k}\\
    &\leq \max(|a_k-\bar{a}|,|b_k-\bar{a}|)|y|+\bar{C}(A_k+B_k)|y|/\eta.
  \end{aligned}
\end{equation*}
Define $\sigma_v(|y|)=\max(|a_k-\bar{a}|,|b_k-\bar{a}|)+\bar{C}(A_k+B_k)/\eta$ and then $\sigma_v(r)\rightarrow 0$ as $r\rightarrow 0$. That is, $v$ (and hence $u$) is differentiable at $0$.~\qed~\\

%
\bibliographystyle{model4-names}
\bibliography{PDE}

\end{document}